\documentclass[11pt,english]{amsart}
\usepackage[T1]{fontenc}
\usepackage[latin9]{inputenc}
\usepackage{amstext}
\usepackage{amsthm}
\usepackage{amssymb}
\usepackage{esint}

\makeatletter
\usepackage{amsmath,amsfonts,amssymb,amsthm,epsfig}

\usepackage{color}
\usepackage[final,allcolors=blue,colorlinks=true]{hyperref}

\def\R{\mathbb R}
\def\N{\mathbb N}

\def\al{\alpha}

\def\ep{\epsilon}
\def\la{\lambda}

\def\var{\varphi}
\def\om{\omega}
\def\na{\nabla}
\def\Om{\Omega}  
\def\De{\Delta}      

\def\cal{\mathcal}
\def\wq{\infty}
\def\pa{\partial}
\def\divergence{\text{\rm div}\,}

\def\loc{\text{\rm loc}}

 \def\diam{\text{\rm diam}}

\newcommand{\medint}{-\kern -,375cm\int}         
\newcommand{\medintinrigo}{-\kern -,315cm\int}

\numberwithin{equation}{section}
\newtheorem{theorem}{Theorem}[section]
\newtheorem*{theorem*}{Theorem}  

\newtheorem*{conclusion*}{Conclusin}

\newtheorem*{corollary*}{Corollary}

\newtheorem*{lemma*}{Lemma}

\newtheorem*{notation*}{Notation}

\newtheorem{proposition}[theorem]{Proposition}
\newtheorem*{proposition*}{Proposition}

\newtheorem*{remark*}{Remark}

\newtheorem*{example*}{Example}                
\theoremstyle{definition}

\makeatother

\usepackage{babel}
\begin{document}
\title[]{Sharp Morrey regularity theory for a fourth order geometrical equation}

\author[C.-L. Xiang]{Chang-Lin Xiang and Gao-Feng Zheng*}

\address{Three Gorges Mathematical Research Center, China Three Gorges University, 443002, Yichang, People's Republic of China}
\email{changlin.xiang@ctgu.edu.cn}

\address[Gao-Feng Zheng]{School of Mathematics and Statistics, Central China Normal University, Wuhan 430079,  P. R.  China}
\email{gfzheng@mail.ccnu.edu.cn}

\thanks{*: corresponding author}
\thanks{The authors are financially supported by the National Natural Science Foundation of China (No. 12271296, 12271195).}

\begin{abstract}
This paper is a continuation of the recent work of Guo-Xiang-Zheng \cite{Guo-Xiang-Zheng-2021-CV}. We deduce sharp Morrey regularity theory for weak solutions to the fourth order nonhomogeneous Lamm-Rivi\`ere equation
\begin{equation*} \De^{2}u=\De(V\na u)+div(w\na u)+(\na\om+F)\cdot\na u+f\qquad\text{in }B^{4},\end{equation*} under  smallest regularity assumptions of $V,w,\om, F$ and that $f$ belongs to some Morrey spaces, which was motivated by many geometrical problems such as the flow of biharmonic mappings. Our results deepens the $L^p$ type regularity theory of  \cite{Guo-Xiang-Zheng-2021-CV}, and generalizes the work of Du, Kang and Wang \cite{Du-Kang-Wang-2022} on  a second order problem to our fourth order problems.
\end{abstract}

\maketitle

{\small
\keywords {\noindent {\bf Keywords:} fourth order elliptic equation, regularity theory, Morrey space, Decay estimates, Riesz potential}
\smallskip
\newline
\subjclass{\noindent {\bf 2020 Mathematics Subject Classification:} 35J47,  35B65}
\tableofcontents}
\bigskip

\section{Introduction and main results}

In this note, we consider the fourth order elliptic problem
\begin{equation}
\De^{2}u=\De(V\na u)+div(w\na u)+(\na\om+F)\cdot\na u+f\qquad\text{in }B^{4},\label{eq: nonhomo-L-R}
\end{equation}
where $B^{4}=B_{1}^{4}(0)\subset\R^{4}$ is the unit ball, $u\in W^{2,2}(B^{4},\R^{m})$
is a weak solution, and $V,w,\om,F$ satisfy the smallest regularity
assumption
\[
\begin{aligned} & V\in W^{1,2}(B^{4},M_{m}\otimes\Lambda^{1}\R^{4}),\quad\quad w\in L^{2}(B^{4},M_{m}),\\
 & \om\in L^{2}(B^{4},so_{m}),\quad\quad F\in L^{\frac{4}{3},1}(B^{4},M_{m}\otimes\Lambda^{1}\R^{4}),
\end{aligned}
\]
which makes the equation  critical when discussing regularity theory.

When $f\equiv0$, the corresponding homogeneous problem
\begin{equation}
\De^{2}u=\De(V\na u)+div(w\na u)+(\na\om+F)\cdot\na u\qquad\text{in }B^{4}\label{eq: homo L-R}
\end{equation}
was first formulated in the interesting work \cite{Lamm-Riviere-2008}
of Lamm and Rivi\`ere, which was intended to extend the celebrated
work of Rivi\`ere \cite{Riviere-2007} on the second order problem
\begin{equation}
-\De u=\Om\cdot\na u\qquad\text{in }B^{2};\label{eq: Riviere}
\end{equation}
and also to provide a new approach to some fourth order conformally
invariant variational problems such as biharmonic mappings (see e.g.
\cite{Chang-W-Y-1999,Wang-2017-CV,Wang-2004-MZ,Wang-2004-CPAM} for
more regularity theory on biharmonic mappings). The reader is also
referred to \cite{Lamm-Riviere-2008} for details on how to write
the equation of biharmonic mappings into the form (\ref{eq: homo L-R}).

The flow of extrinsic or intrinsic biharmonic mappings also attracted
extensive study, see e.g. Gastel \cite{Gastel-2006-AG}, Lamm \cite{Lamm-2004-AGAG},
Wang \cite{Wang-2007}. In this case, one is naturally motivated to
consider the nonhomogeneous problem (\ref{eq: nonhomo-L-R}). For
instance, using the advantage of equation (\ref{eq: nonhomo-L-R}),
the existence of global flow of extrinsic biharmonic mappings was
obtained in \cite{Lamm-Riviere-2008} for any initial data. We also
mention that another slightly different fourth order elliptic system
but with the same nature as (\ref{eq: homo L-R}) was formulated by
Struwe \cite{Struwe-2008} and further studied by Guo, Wang and Xiang \cite{Guo-Wang-Xiang-2022-CV} and many others. The main results of Lamm and Rivi\`ere
\cite{Lamm-Riviere-2008} can be summarized as follows.

\begin{theorem} \label{thm: Lamm-Riviere} (\cite{Lamm-Riviere-2008})
(1) For any $m\in\N$, there exist  constants $C_{m}>0$ and $\ep_{m}>0$
 such that  if
\begin{equation}
\|V\|_{W^{1,2}(B^{4})}+\|w\|_{L^{2}(B^{4})}+\|\om\|_{L^{2}(B^{4})}+\|F\|_{L^{4/3,1}(B^{4})}<\ep_{m},\label{eq: smallness assumption}
\end{equation}
then there exist  $A\in W^{2,2}\cap L^{\wq}(B_{1/2}^{4},M(m))$ and
$B\in W^{1,4/3}(B_{1/2}^{4},M(m)\otimes\wedge^{2}\R^{4})$ with
\begin{equation}
\|A\|_{W^{2,2}(B_{1/2}^{4})}+\|{\rm dist}(A,SO_{m})\|_{L^{\wq}(B_{1/2}^{4})}+\|B\|_{W^{1,4/3}(B_{1/2}^{4})}\le C_{m}\ep_{m}\label{eq: A-B small}
\end{equation}
such that
\[
\na\De A+\De AV-\na Aw+A(\na\om+F)=\sum_{l,k=1}^{4}\pa_{l}B_{lk}\pa_{k}.
\]

(2) Consequently, under the smallness assumption (\ref{eq: smallness assumption}),
$u$ solves problem (\ref{eq: homo L-R}) in $B_{1/2}^{4}$ if and
only if $u$ satisfies the conservation law
\begin{equation}
\De(A\De u)=\divergence(K)\qquad\text{in }B_{1/2}^{4},\label{eq:conservation law of Lamm Riviere}
\end{equation}
where
\begin{equation}
K\equiv2\na A\De u-\De A\na u+Aw\na u-\na AV\cdot\na u+A\na(V\cdot\na u)+B\cdot\na u.\label{eq: K}
\end{equation}

(3) Every weak solution of equation (\ref{eq: homo L-R}) is continuous.
\end{theorem}

For applications of this result to the flow of existence biharmonic
mappings, the reader is referred to the same work \cite{Lamm-Riviere-2008}.

As an improvement, Guo and Xiang \cite{Guo-Xiang-2019-Boundary,Guo-Xiang-2019-Higher}
showed that every weak solution of equation (\ref{eq: homo L-R})
is in fact locally H\"older continuous. Quite recently, Guo, Xiang
and Zheng \cite{Guo-Xiang-Zheng-2021-CV} considered $L^{p}$ type
regularity theory of problem (\ref{eq: nonhomo-L-R}). They proved
that if $f\in L^{p}(B^{4})$ for some $1<p<4/3$, then every weak
solution $u\in W_{\loc}^{3,p_{1}}(B^{4})$, where $p_{1}=4p/(4-p)$;
while in the case $p=1$ and $f\in L\log L(B^{4})$ there holds $u\in W_{\loc}^{3,4/3,1}(B^{4})$.
The result is optimal in the sense that, the regularity of $u$ can
not be improved too much even if $f\in L^{p}(B^{4})$ for some $p\ge4/3$.
Also, the third order derivative of $u$ is the best possible: counterexample
shows that in general weak solutions may not have fourth order weak
derivatives.

As a counterpart of the second order problem (\ref{eq: Riviere}),
we mention that Rivi\`ere \cite{Riviere-2007} was the first one
to formulate problem (\ref{eq: Riviere}) in order to study second
order conformally invariant variational problems, and he also established
the method of conservation law and continuity of weak solutions to
problem (\ref{eq: Riviere}). Sharp and Topping \cite{Sharp-Topping-2013-TAMS}
considered the nonhomogeneous problem
\begin{equation}
-\De u=\Om\cdot\na u+f\qquad\text{in }B^{2}.\label{eq: nonhomo-R}
\end{equation}
They proved, among other results, that if $f\in L^{p}(B^{2})$ for
some $1<p<2$, then every weak solution $u\in W^{1,2}(B^{2})$ of
problem (\ref{eq: nonhomo-R}) belongs to $W_{\loc}^{2,p}(B^{2})$;
and in the case $p=1$ and $f\in L\log L(B^{2})$, then $u\in W_{\loc}^{2,1}(B^{2})$.
Inspired by the work of Wang \cite{Wang-2017-CV} on harmonic mappings,
recently Du, Kang and Wang \cite{Du-Kang-Wang-2022} extend the theory
of Sharp and Topping \cite{Sharp-Topping-2013-TAMS} to the case that
$f$ belongs to some Morrey spaces (for definitions of Morrey spaces,
see Appendix \ref{sec: Morrey-spaces}). They proved that if $f\in M^{1,\la}(B^{2})$
for some $0<\la<1$, then $\na u\in M_{\loc}^{2,2\la}(B^{2})$; and
if $f\in M^{p,p\la}(B^{2})$ for some $1<p<2$ and $0\le\la<(2-p)/p$,
then $\na^{2}u\in M_{\loc}^{p,p\la}(B^{2})$ and $\na u\in M_{\loc}^{p^{*},p^{*}\la}(B^{2})$.
This also extends the H\"older regularity theory of Wang \cite{Wang-2017-CV}
on harmonic mappings to the general second order problem (\ref{eq: nonhomo-R}),
which maybe have potential applications to nonhomogeneous prescribed
mean curvature equations.

This note can be viewed as a continuation of \cite{Guo-Xiang-Zheng-2021-CV}.
Our motivation is to further extend the $L^{p}$ type regularity theory
of \cite{Guo-Xiang-Zheng-2021-CV} to the case that $f$ belongs to
some Morrey spaces, and thus extend the Morrey regularity theory of
Du, Kang and Wang \cite{Du-Kang-Wang-2022} to the fourth order problem
(\ref{eq: nonhomo-L-R}). Out first result reads as follows.

\begin{theorem}\label{thm: Morrey decay} Suppose $u\in W^{2,2}(B^{4})$
is a weak solution to equation (\ref{eq: nonhomo-L-R}) with $f\in M^{1,\la}(B^{4})$
for some $0<\la<1$. Then $\na u\in M_{\loc}^{4,4\la}(B^{4})$, $\na^{2}u\in M_{\loc}^{2,2\la}(B^{4})$
with
\begin{equation}
\|\na u\|_{M^{4,4\la}(B_{1/2}^{4})}+\|\na^{2}u\|_{M^{2,2\la}(B_{1/2}^{4})}\le C\left(\|u\|_{W^{2,2}(B^{4})}+\|f\|_{M^{1,\la}(B^{4})}\right).\label{eq: Morrey estimate 1}
\end{equation}
As a result, we have $u\in C_{\loc}^{0,\la}(B^{4})$ with
\begin{equation}
\|u\|_{C^{0,\la}(B_{1/2}^{4})}\le C\left(\|u\|_{W^{2,2}(B^{4})}+\|f\|_{M^{1,\la}(B^{4})}\right),\label{eq: Holder continuity 1}
\end{equation}
where $C>0$ is a positive constant depending only on $m,\la$ and
the coefficient functions $V,w,\om,F$.\end{theorem}

For the definition of Morrey spaces, see Appendix \ref{sec: Morrey-spaces}.
The result is optimal in the following sense. Consider the simplest
case $\De^{2}u=f\in L^{p}(B^{4})$ for some $1<p<4/3$. Then $f\in L^{p}(B^{4})\subset M^{1,\la}(B^{4})$
with $\la=4(1-1/p)\in(0,1)$. By the standard elliptic regularity
theory, $u\in W_{\loc}^{4,p}(B^{4})\subset C_{\loc}^{0,\la}(B^{4})$.
This shows that the H\"older regularity (\ref{eq: Holder continuity 1})
is optimal. So it is natural to ask whether this result hold for $\la=1$.
The answer is negative, even under stronger assumptions on $f$. Indeed,
consider the case $p=4/3$. Then $f\in L^{4/3}(B^{4})\subset M^{1,1}(B^{4})$.
It is known from elliptic regularity theory that the best possible
regularity is $u\in W_{\loc}^{4,4/3}(B^{4})$ and it is possible that
$u\not\in C_{\loc}^{0,1}(B^{4})$.

In the case $f$ has better Morrey regularity, we have the following
improved regularity.

\begin{theorem} \label{thm: Sobolev-Morrey-regularity} Suppose $u\in W^{2,2}(B^{4})$
is a weak solution to equation (\ref{eq: nonhomo-L-R}) with $f\in M^{p,p\la}(B_{1}^{4})$
for some $1<p<4/3$ and $0\le\la<(4-3p)/p$. Then $\na u\in M_{\loc}^{p_{3},p_{3}\la}(B^{4})$,
$\na^{2}u\in M_{\loc}^{p_{2},p_{2}\la}(B^{4})$, $\na^{3}u\in M_{\loc}^{p_{1},p_{1}\la}(B^{4})$
with $p_{i}=4p/(4-ip)$ ($1\le i\le3$), and
\[
\sum_{i=1}^{3}\left\Vert \na^{i}u\right\Vert _{M^{p_{4-i},p_{4-i}\la}(B_{1/2}^{4})}\le C\left(\|u\|_{W^{2,2}(B^{4})}+\|f\|_{M^{1,\la}(B^{4})}\right)
\]
holds for some constant $C>0$ depending only on $m,p,\la$ and the
coefficient functions $V,w,\om,F$. \end{theorem}

As $M^{p,0}(B^{4})=L^{p}(B^{4})$, this result can be viewed as an
extension of the $L^{p}$ type regularity theory of Guo, Xiang and
Zheng \cite{Guo-Xiang-Zheng-2021-CV}. As already explained in \cite{Guo-Xiang-Zheng-2021-CV},
in general, the third order Sobolev regularity is the best possible.
Thus in our case we still have as most third order regularities. Moreover,
the indexes are all optimal, which can be seen from the simple case
$\De^{2}u=f$. We leave the details for interested readers.

The strategy of proving Theorems \ref{thm: Morrey decay} and \ref{thm: Sobolev-Morrey-regularity}
can be briefly illustrated as follows. Theorems \ref{thm: Morrey decay}
follows from a Morrey type decay estimates, in which the conservation
law of Theorem \ref{thm: Lamm-Riviere} is crucial. To prove Theorem
\ref{thm: Sobolev-Morrey-regularity}, we will begin with the $W^{3,p_{1}}$-regularity
theory of (\ref{eq: nonhomo-L-R}) in view of \cite{Guo-Xiang-Zheng-2021-CV}.
So under the assumption $f\in M^{p,p\la}(B^{4})$, we have $u\in W^{3,p_{1}}(B_{3/4}^{4})$.
Next we use a perturbation method and duality method to derive the
Morrey regularity of $\na^{i}u$ ($i=1,2,3$).

Our notations are standard. Throughout we use $A\lesssim B$ to mean
that $A\le CB$ for some constant $C>0$ depending only $m,p,\la$
and independent of the parameter $\ep$. We use $I_{\al}=|x|^{\al-n}$
to denote Riesz potential for all $0<\al<n$, such that by $I_{\al}(f)$
we mean
\[
I_{\al}(f)=|x|^{\al-n}\ast f.
\]
For Riesz potential theory on Morrey spaces, see Proposition (\ref{prop: Adams})
in the Appendix.

\section{Proof of theorem \ref{thm: Morrey decay} }

In this section we prove theorem \ref{thm: Morrey decay} using a
Morrey type decay estimates.

Using a standard translation and scaling argument (see \cite[subsection 2.3]{Guo-Xiang-Zheng-2021-CV}
for details), we can assume
\[
\|V\|_{W^{1,2}(B^{4})}+\|w\|_{L^{2}(B^{4})}+\|\om\|_{L^{2}(B^{4})}+\|F\|_{L^{4/3,1}(B^{4})}<\ep,
\]
where $\ep<\ep_{m}$ so that by the conservation law of Theorem \ref{thm: Lamm-Riviere}
we have
\[
\De{\rm div}(A\na u)={\rm div}\tilde{K}+Af\qquad\text{in }B_{1/2}^{4},
\]
where
\[
\tilde{K}=K+\De(\na A\cdot\na u)
\]
with $A,B$ and $K$ being defined as in Theorem \ref{thm: Lamm-Riviere}.

To study $A\na u$, we use the Hodge decomposition to deduce
\[
Adu=d\tilde{u}_{1}+d^{\ast}\tilde{u}_{2}+\tilde{h}\qquad\text{in }B_{1/2}^{4}
\]
for some harmonic 1-form $\tilde{h}$, where $\tilde{u}_{1}$ is a
function and $\tilde{u}_{2}$ is a $2$-form. Then $\tilde{u}_{1},\tilde{u}_{2}$
satisfy equations
\[
\De^{2}\tilde{u}_{1}={\rm div}\tilde{K}+Af\qquad\text{in }B_{1/2}^{4},
\]
\[
\De\tilde{u}_{2}=dA\wedge du\qquad\text{in }B_{1/2}^{4}.
\]

To proceed, let us extend $V,w,\om,F,A,B$ from $B_{1/2}^{4}$ to
the whole space $\R^{4}$ such that for some $C_{m}>0$ independent
of $\ep$, there hold
\[
\|V\|_{W^{1,2}(\R^{4})}+\|w\|_{L^{2}(\R^{4})}+\|\om\|_{L^{2}(\R^{4})}+\|F\|_{L^{4/3,1}(\R^{4})}<C_{m}\ep,
\]
\begin{equation}
\|A\|_{W^{2,2}(\R^{4})}+\|B\|_{W^{1,4/3}(\R^{4})}<C_{m}\ep\label{eq: estimate of A-B in R4}
\end{equation}
(still denoted by the same notations). As to $f$, we simply let $f\equiv0$
in $\R^{4}\backslash B_{1/2}^{4}$ so that
\[
\|f\|_{M^{1,\la}(\R^{4})}\le C(\la)\|f\|_{M^{1,\la}(B_{1/2}^{4})}
\]
 for some constant $C(\la)$ depending only on $\la$.

Put
\[
u_{11}=c\log\ast\left({\rm div}\tilde{K}\right),\qquad u_{12}=c\log\ast(Af),\quad\quad u_{2}=I_{2}(dA\wedge du),
\]
where $c\log$ and $I_{2}=c|x|^{-2}$ denote the fundamental solutions
of $\De^{2}$ and $\De$ in $\R^{4}$ respectively, so that
\[
\De^{2}u_{11}={\rm div}\tilde{K},\quad\quad\De^{2}u_{12}=Af,\quad\quad\De u_{2}=dA\wedge du
\]
in $\R^{4}$. Then we find that $\De^{2}(u_{11}+u_{12}-u_{1})=0$
and $\De(u_{2}-\tilde{u}_{2})=0$ in $B_{1/2}^{4}$ respectively.
Thus, we obtain the decomposition
\begin{equation}
Adu=du_{11}+du_{12}+d^{\ast}u_{2}+h\label{eq: Hodge}
\end{equation}
for some biharmonic function $h$ in $B_{1/2}^{4}$.

Now we can use Riesz potential theory to estimate each term in equation
(\ref{eq: Hodge}).

For the first term $u_{11}$, note that
\[
|\na u_{11}|\lesssim I_{2}(|\tilde{K}|).
\]
$\tilde{K}$ contains 7 terms. For the first term $\na A\De u$, we
have
\[
\left\Vert I_{2}(\na A\De u)\right\Vert _{L^{4}(\R^{4})}\lesssim\|\na A\|_{L^{4}(\R^{4})}\|\na^{2}u\|_{L^{2}(\R^{4})}\lesssim\ep\|\na^{2}u\|_{L^{2}(\R^{4})},
\]
where we used the boundedness of the operator $I_{2}:L^{4/3}(\R^{4})\to L^{4}(\R^{4})$
and the estimate (\ref{eq: estimate of A-B in R4}). Similarly we
can estimate the remaining six terms of $\na u_{11}$ to find that
\[
\left\Vert \na u_{11}\right\Vert _{L^{4}(\R^{4})}\lesssim\ep\left(\|\na^{2}u\|_{L^{2}(\R^{4})}+\|\na u\|_{L^{4}(\R^{4})}\right).
\]
Since $u$ is extended to $\R^{4}$ in a bounded way from $B_{1/2}^{4}$,
we obtain
\begin{equation}
\left\Vert \na u_{11}\right\Vert _{L^{4}(\R^{4})}\lesssim\ep\left(\|\na^{2}u\|_{L^{2}(B_{1/2}^{4})}+\|\na u\|_{L^{4}(B_{1/2}^{4})}\right).\label{eq: estimate of u11}
\end{equation}

To estimate $u_{12}$, we use the estimate
\[
|\na u_{12}|\lesssim I_{3}(|\na A||\na u|)
\]
and the boundedness of the operator (see Proposition \ref{prop: Adams}
in Appendix \ref{sec: Morrey-spaces})
\[
I_{3}:M^{1,\la}(\R^{4})\to M_{\ast}^{\frac{4-\la}{1-\la},\la}(\R^{4})
\]
and find that $\na u_{12}\in M_{\ast}^{\frac{4-\la}{1-\la},\la}(\R^{4})$
with estimate
\[
\left\Vert \na u_{12}\right\Vert _{M_{\ast}^{\frac{4-\la}{1-\la},\la}(\R^{4})}\lesssim\|f\|_{M^{1,\la}(\R^{4})}\le C(\la)\|f\|_{M^{1,\la}(B_{1/2}^{4})}.
\]
In particular, for any $r>0$, we obtain by H\"older's inequality
that
\[
\|\na u_{12}\|_{L^{4}(B_{r}^{4})}\le C\left\Vert u_{12}\right\Vert _{L_{\ast}^{\frac{4-\la}{1-\la}}(B_{r}^{4})}r^{4\left(\frac 14-\frac{1-\la}{4-\la}\right)}\le C\|f\|_{M^{1,\la}(B_{1/2}^{4})}r^{\la}.
\]
Thus, $u_{12}\in M^{4,4\la}(\R^{4})$ with
\begin{equation}
\left\Vert \na u_{12}\right\Vert _{M^{4,4\la}(\R^{4})}\le C\|f\|_{M^{1,\la}(B_{1/2}^{4})}.\label{eq: estimate of u12}
\end{equation}

The estimate of $u_{2}$ is totally the same as above, which gives
\begin{equation}
\left\Vert \na u_{2}\right\Vert _{L^{4}(\R^{4})}\lesssim\ep\|\na u\|_{L^{4}(\R^{4})}\lesssim\ep\left(\|\na^{2}u\|_{L^{2}(B_{1/2}^{4})}+\|\na u\|_{L^{4}(B_{1/2}^{4})}\right).\label{eq: estimate of u2}
\end{equation}

Since $h$ is biharmonic in $B_{1/2}^{4}$, for any $1/4>r>0$, we
have
\begin{equation}
\|h\|_{L^{4}(B_{r}^{4})}\le Cr\|h\|_{L^{4}(B_{1/2}^{4})}.\label{eq: estimate of h}
\end{equation}

Therefore, combining the above estimates (\ref{eq: estimate of u11}),
(\ref{eq: estimate of u12}), (\ref{eq: estimate of u2}) and (\ref{eq: estimate of h})
on $u_{11}$, $u_{12}$, $u_{2}$ and $h$, we obtain, for any $1/4>r>0$,
\[
\begin{aligned}\|\na u\|_{L^{4}(B_{r}^{4})} & \le\|h\|_{L^{4}(B_{r}^{4})}+\|\na u_{11}\|_{L^{4}(B_{r}^{4})}+\|\na u_{12}\|_{L^{4}(B_{r}^{4})}+\|\na u_{2}\|_{L^{4}(B_{r}^{4})}\\
 & \lesssim r\|h\|_{L^{4}(B_{1/2}^{4})}+\ep\left(\|\na^{2}u\|_{L^{2}(B_{1/2}^{4})}+\|\na u\|_{L^{4}(B_{1/2}^{4})}\right)+\|f\|_{M^{1,\la}(B_{1/2}^{4})}r^{\la}\\
 & \le C(r+\ep)\left(\|\na^{2}u\|_{L^{2}(B_{1/2}^{4})}+\|\na u\|_{L^{4}(B_{1/2}^{4})}\right)+C\|f\|_{M^{1,\la}(B_{1/2}^{4})}r^{\la},
\end{aligned}
\]
where $C$ is a constant depending only on $m,\la$.

Similarly, repeating the above procedure, we can obtain the same estimate
for $\|\na^{2}u\|_{L^{2}(B_{r}^{4})}$. Hence, by writing
\[
U(r)\equiv\|\na u\|_{L^{4}(B_{r}^{4})}+\|\na^{2}u\|_{L^{2}(B_{r}^{4})},
\]
we derive for any $0<r<1/4$
\[
U(r)\le C(r+\ep)U(1/2)+C\|f\|_{M^{1,\la}(B_{1/2}^{4})}r^{\la}
\]
for some constant $C>0$ depending only on $m$ and $\la$. We first
$r_{0}\ll1/4$ such that $2Cr_{0}\le r_{0}^{(\la+1)/2}$, and then
choose $\ep\le r_{0}$ so that
\[
U(r_{0})\le U(1/2)r_{0}^{(\la+1)/2}+C\|f\|_{M^{1,\la}(B_{1/2}^{4})}r_{0}^{\la}.
\]
 Now using a standard scaling and iteration argument (see \cite[Theorem 3.1]{Guo-Xiang-Zheng-2021-CV}
for instance), we obtain, for any $k\ge1$,
\[
U(r_{0}^{k})\le r_{0}^{\frac{\la+1}{2}}U(r_{0}^{k-1})+C\|f\|_{M^{1,\la}(B_{1/2}^{4})}r_{0}^{k\la},
\]
which implies that
\[
U(r)\le C\left(U(1)+\|f\|_{M^{1,\la}(B_{1}^{4})}\right)r^{\la}
\]
for all $0<r<1/4$. This proves (\ref{eq: Morrey estimate 1}).

The H\"older continuity estimate (\ref{eq: Holder continuity 1})
follows from a standard Morrey type Dirichlet growth theorem, see
Giaquinta \cite{Giaquinta-Book} for example. This completes the proof
of Theorem \ref{thm: Morrey decay}.

\section{Proof of theorem \ref{thm: Sobolev-Morrey-regularity}}

In this section we prove Theorem \ref{thm: Sobolev-Morrey-regularity}.

Suppose now $u\in W^{2,2}(B^{4})$ is a weak solution to equation
(\ref{eq: nonhomo-L-R}) with $f\in M^{p,p\la}(B_{1}^{4})$ for some
$1<p<4/3$ and $0<\la<(4-3p)/p$. Since the result is local and the
problem is invariant under suitable scaling (for details, see \cite{Guo-Xiang-Zheng-2021-CV}),
we can assume that
\[
\|V\|_{W^{1,2}(B_{1})}+\|w\|_{L^{2}(B_{1})}+\|\om\|_{L^{2}(B_{1})}+\|F\|_{L^{4/3,1}(B_{1})}<\ep_{m}
\]
for some $\ep_{m}\ll1$. Under this assumption, let us first improve
the result of Theorem \ref{thm: Morrey decay} for later use.

1. By H\"older's inequality, we easily find that
\begin{equation}
M^{p,p\la}(B_{1}^{4})\subset M^{1,\la_{0}}(B_{1}^{4}),\label{eq: Morrey-Holder embedding}
\end{equation}
where
\begin{equation}
\la_{0}=\la+4(1-1/p).\label{eq: lambda-0}
\end{equation}
Thus by Theorem \ref{thm: Morrey decay}, we have
\[
\|\na u\|_{M^{4,4\la_{0}}(B_{1/2}^{4})}+\|\na^{2}u\|_{M^{2,2\la_{0}}(B_{1/2}^{4})}\le C\left(\|u\|_{W^{2,2}(B^{4})}+\|f\|_{M^{1,\la_{0}}(B^{4})}\right)
\]
and $u\in C^{0,\la_{0}}(B_{1/2}^{4})$ with
\begin{equation}
\|u\|_{C^{0,\la_{0}}(B_{1/2}^{4})}\le C\left(\|u\|_{W^{2,2}(B^{4})}+\|f\|_{M^{1,\la_{0}}(B^{4})}\right).\label{eq: Holder-2}
\end{equation}

2. Since $f\in L^{p}(B^{4})$ for $1<p<4/3$, by Guo-Xiang-Zheng \cite{Guo-Xiang-Zheng-2021-CV},
$u\in W^{3,p_{1}}(B_{1/2}^{4})$ with
\[
\|u\|_{W^{3,p_{1}}(B_{1/2}^{4})}\le C\left(\|u\|_{W^{2,2}(B^{4})}+\|f\|_{L^{p}(B^{4})}\right),
\]
where
\[
p_{1}=4p/(4-p).
\]

3. Due to the above $W^{3,p_{1}}$-regularity, we can then repeat
the proof of Theorem \ref{thm: Morrey decay} to find that
\[
\|\na^{3}u\|_{M^{\frac{4}{3},\frac{4}{3}\la_{0}}(B_{1/2}^{4})}\le C\left(\|u\|_{W^{2,2}(B^{4})}+\|f\|_{M^{1,\la_{0}}(B^{4})}\right).
\]
Hence, we summarize the above result to conclude that
\begin{equation}
\|\na u\|_{M^{4,4\la_{0}}(B_{1/2}^{4})}+\|\na^{2}u\|_{M^{2,2\la_{0}}(B_{1/2}^{4})}+\|\na^{3}u\|_{M^{\frac{4}{3},\frac{4}{3}\la_{0}}(B_{1/2}^{4})}\le CM,\label{eq: Morrey decay-2}
\end{equation}
where
\begin{equation}
M=\|u\|_{W^{2,2}(B^{4})}+\|f\|_{M^{p,p\la}(B^{4})}.\label{eq: M}
\end{equation}

Now we can prove Theorem \ref{thm: Sobolev-Morrey-regularity}. Due
to the second remark above, we know now that $u\in W^{3,p_{1}}(B_{1/2}^{4})$.
Our aim is to deduce its Morrey regularity.

Fix $x_{0}\in B_{1/4}^{4}$ and $0<r<1/4$. Split $A\De u=v+h$ in
$B_{r}(x_{0})$ with $v$ and $h$ satisfying
\begin{equation}
\begin{cases}
\De v={\rm div}K+Af, & \text{in }B_{r}(x_{0}),\\
v=0 & \text{on }\pa B_{r}(x_{0}),
\end{cases}\label{eq: v}
\end{equation}
and
\begin{equation}
\begin{cases}
\De h=0 & \text{in }B_{r}(x_{0}),\\
h=A\De u & \text{on }\pa B_{r}(x_{0}).
\end{cases}\label{eq: h}
\end{equation}

1. Estimate of $v$.

Write $p_{2}=4p/(4-2p)$ and $p_{3}=4p/(4-3p)$. By the duality of
$L^{p}$-space, we have

\[
\|v\|_{L^{p_{2}}(B_{r}(x_{0}))}=\sup_{\var\in{\cal A}}\int_{B_{r}(x_{0})}v\var,
\]
where
\[
{\cal A}=\{\var\in L^{p_{2}^{\prime}}(B_{r}(x_{0});\R^{m}):\|\var\|_{p_{2}^{\prime}}\le1\}
\]
and $p_{2}^{\prime}=p_{2}/(p_{2}-1)$ is the H\"older conjugate of
$p_{2}$. For any $\var\in{\cal A}$, the standard elliptic regularity
theory implies that there is a unique $\psi\in W^{2,p_{2}^{\prime}}\cap W_{0}^{1,p_{2}^{\prime}}(B_{r}(x_{0}))$
satisfying $\De\psi=\var$ in $B_{r}(x_{0})$, $\psi=0$ on $\pa B_{r}(x_{0})$,
and moreover,
\[
\|\psi\|_{W^{2,p_{2}^{\prime}}(B_{r}(x_{0}))}\le C(p).
\]
The Sobolev embedding theorem implies that
\[
W^{2,p_{2}^{\prime}}(B_{r}(x_{0}))\subset W^{1,\frac{4p_{2}^{\prime}}{4-p_{2}^{\prime}}}(B_{r}(x_{0}))\subset L^{p^{\prime}}(B_{r}(x_{0})).
\]
Hence
\[
\left\Vert \psi\right\Vert _{L^{p^{\prime}}(B_{r}(x_{0}))}+\left\Vert \na\psi\right\Vert _{L^{4p_{2}^{\prime}/(4-p_{2}^{\prime})}(B_{r}(x_{0}))}\le C(p).
\]

Now integrating by parts yields
\[
\int_{B_{r}(x_{0})}v\var=\int_{B_{r}(x_{0})}v\De\psi=-\int_{B_{r}(x_{0})}K\na\psi+\int_{B_{r}(x_{0})}Af\psi=:I_{1}+I_{2}.
\]
Using H\"older's inequality we easily obtain
\[
|I_{2}|\le C\|f\|_{L^{p}(B_{r}(x_{0}))}\|\psi\|_{L^{p^{\prime}}(B_{r}(x_{0}))}\le C\|f\|_{M^{p,p\la}(B^{4})}r^{\la}.
\]
Here we used the assumption that $f\in M^{p,p\la}(B^{4})$ in the
last inequality. To estimate $I_{1}$, we use H\"older's inequality
again to infer that
\[
\int_{B_{r}(x_{0})}|\na A\De u\na\psi|\le\|\na A\|_{4}\|\|\De u\|_{L^{p_{2}}(B_{r}(x_{0}))}\left\Vert \na\psi\right\Vert _{L^{4p_{2}^{\prime}/(4-p_{2}^{\prime})}(B_{r}(x_{0}))}\le C\ep\|\De u\|_{L^{p_{2}}(B_{r}(x_{0}))}.
\]
Other terms of $K$ can be estimated similarly. This gives
\[
|I_{1}|\le C\ep\left(\|\De u\|_{L^{p_{2}}(B_{r}(x_{0}))}+\|\na u\|_{L^{p_{3}}(B_{r}(x_{0}))}\right)
\]
Hence combining the above two estimates on $I_{1}$ and $I_{2}$ gives
\[
\|v\|_{L^{p_{2}}(B_{r}(x_{0}))}\le C\ep\left(\|\De u\|_{L^{p_{2}}(B_{r}(x_{0}))}+\|\na u\|_{L^{p_{3}}(B_{r}(x_{0}))}\right)+C\|f\|_{M^{p,p\la}(B^{4})}r^{\la}.
\]

We can estimate the norm of $\na v$ similarly. Note that
\begin{equation}
\|\na v\|_{L^{p_{1}}(B_{r}(x_{0}))}=\sup_{\var\in{\cal B}}\int_{B_{r}(x_{0})}\na v\cdot\var,\label{eq: duality of nabla v}
\end{equation}
where
\[
{\cal B}=\left\{ \var\in L^{p_{1}^{\prime}}(B_{r}(x_{0}),\R^{m}\otimes\R^{4}):\|\var\|_{L^{p_{1}^{\prime}}(B_{r}(x_{0}))}\le1\right\} .
\]
Using Helmholtz decomposition, for each $\var\in{\cal B}$, there
exist $\psi\in W_{0}^{1,p_{1}^{\prime}}(B_{r}(x_{0}),\R^{m})$ and
$\xi\in L^{p_{1}^{\prime}}(B_{r}(x_{0}),\R^{m}\otimes\R^{4})$ such
that $\var=\na\psi+\xi$ and ${\rm div}\xi=0$ in $B_{r}(x_{0})$.
Substituting this decomposition into (\ref{eq: duality of nabla v})
gives
\[
\int_{B_{r}(x_{0})}\na v\cdot\var=\int_{B_{r}(x_{0})}\na v\cdot\na\psi
\]
in view of the boundary condition $v=0$ on $\pa B_{r}(x_{0})$. Working
as in the above for $\|v\|_{L^{p_{2}}(B_{r}(x_{0}))}$, we deduce
\[
\|\na v\|_{L^{p_{1}}(B_{r}(x_{0}))}\le C\ep\left(\|\De u\|_{L^{p_{2}}(B_{r}(x_{0}))}+\|\na u\|_{L^{p_{3}}(B_{r}(x_{0}))}\right)+C\|f\|_{M^{p,p\la}(B^{4})}r^{\la}.
\]

Therefore, we conclude that
\begin{equation}
\begin{aligned} & \qquad\|\na v\|_{L^{p_{1}}(B_{r}(x_{0}))}+\|v\|_{L^{p_{2}}(B_{r}(x_{0}))}\\
 & \le C\ep\left(\|\De u\|_{L^{p_{2}}(B_{r}(x_{0}))}+\|\na u\|_{L^{p_{3}}(B_{r}(x_{0}))}\right)+C\|f\|_{M^{p,p\la}(B^{4})}r^{\la}.
\end{aligned}
\label{eq: estimate of v}
\end{equation}

2. Estimate of $h$.

Since $A\De u\in W^{1,4/3}(B_{r}(x_{0}))$, using a standard result
(see e.g. Lemma C.1 of \cite{Guo-Xiang-Zheng-2021-CV}) we obtain
that $h\in W^{1,4/3}(B_{r}(x_{0}))$ with
\[
\|h\|_{W^{1,4/3}(B_{r}(x_{0}))}\le C\|A\De u\|_{W^{1,4/3}(B_{r}(x_{0}))}\le C\Phi(r)r^{\la_{0}},
\]
where $\la_{0}$ are defined as in (\ref{eq: lambda-0}) and
\[
\Phi(r)=\|\na\De u\|_{L^{p_{1}}(B_{r}(x_{0}))}+\|\De u\|_{L^{p_{2}}(B_{r}(x_{0}))}+\|\na u\|_{L^{p_{3}}(B_{r}(x_{0}))}.
\]
In particular, this implies that
\[
\|h\|_{W^{1,4/3}(B_{r}(x_{0}))}\le CMr^{\la_{0}},
\]
where $M$ is defined as in (\ref{eq: M}). Since $h$ is harmonic
in $B_{r}(x_{0})$, this implies the interior growth estimate
\[
\|h\|_{L^{p_{2}}(B_{r/2}(x_{0}))}\le Cr^{\frac{4}{p_{2}}}\|h\|_{L^{\wq}(B_{r/2}(x_{0}))}\le Cr^{\frac{4}{p_{2}}}\left(\fint_{B_{r}(x_{0})}|h|^{2}\right)^{\frac{1}{2}}\le CMr^{\la}
\]
and similarly
\[
\|\na h\|_{L^{p_{1}}(B_{r/2}(x_{0}))}\le CMr^{\la}.
\]

3. before conclusion, we need to estimate $\na u$. Applying the Sobolev
embedding $W^{2,p_{2}}(B_{r/2}^{4})\subset W^{1,p_{3}}(B_{r/2}^{4})$
we find that
\[
\begin{aligned}\|\na u\|_{L^{p_{3}}(B_{r/4}(x_{0}))} & \le C\|\na^{2}u\|_{L^{p_{2}}(B_{r/4}(x_{0}))}+Cr^{-2}\|u-u_{B_{r/4}(x_{0})}\|_{L^{p_{2}}(B_{r/4}(x_{0}))}\\
 & \le C\|\De u\|_{L^{p_{2}}(B_{r/2}(x_{0}))}+Cr^{-4(1-1/p)}\|u-u_{B_{r/2}(x_{0})}\|_{L^{\wq}(B_{r/2}(x_{0}))}
\end{aligned}
\]
for some constant $C>0$ independent of $r$. Then use the H\"older
estimate (\ref{eq: Holder-2}) we get
\begin{equation}
\|\na u\|_{L^{p_{3}}(B_{r/4}(x_{0}))}\le C\left(\|v\|_{L^{p_{2}}(B_{r/2}(x_{0}))}+\|h\|_{L^{p_{2}}(B_{r/2}(x_{0}))}\right)+CMr^{\la}.\label{eq: estimate of nabla u}
\end{equation}

Finally, by a slight refinement of the above argument, combining the
estimates (\ref{eq: estimate of v}) of $v$ and $h$ and the estimate
(\ref{eq: estimate of nabla u}) of $\na u$, we conclude that
\[
\Phi(r/2)\le C\ep\Phi(r)+CMr^{\la}.
\]
At this moment we use Simon's iteration lemma (see e.g. Lemma A.7
of \cite{Sharp-Topping-2013-TAMS}) to infer that there are $\ep_{0}$
and $r_{0}(\la,p)>0$ sufficiently small such that for all $r\le r_{0}$,
we have
\[
\Phi(r)\le CMr^{\la}.
\]
This completes the proof.

\appendix

\section{Morrey spaces and Riesz potential\label{sec: Morrey-spaces} }

Let $\Om\subset\R^{n}$ be a domain and let $1\le p<\wq$ and $0\le s\le n$.
Denote by $B_{r}(x)$ the ball in $\R^{n}$ centered at $x\in\R^{n}$
with radius $r$. The Morrey space $M^{p,s}(\Om)$ consists of functions
$f\in L^{p}(\Om)$ such that
\[
\|f\|_{M^{p,s}(\Om)}\equiv\sup_{x\in\Om,0<r<\diam\Om}r^{-s/p}\|f\|_{L^{p}(B_{r}(x)\cap\Om)}<\wq.
\]
The weak Morrey space $M_{\ast}^{p,s}(\Om)$ consists of weakly $p$-integrable
functions $f\in L_{\ast}^{p}(\Om)$ such that
\[
\|f\|_{M_{\ast}^{p,s}(\Om)}\equiv\sup_{x\in D,0<r<\diam\Om}r^{-s/p}\|f\|_{L_{\ast}^{p}(B_{r}(x)\cap\Om)}<\wq.
\]
Up to a constant multiplier, the upper bound of $r$ in the supremum
above can be also replaced by $(\diam\Om)/2$ for instance. For more
details, see e.g. \cite{Adams-book}. The following potential theory
is due to Adams \cite{Adams-1975}.

\begin{proposition}\label{prop: Adams} Let $0<\al<n$ and $I_{\al}=|x|^{\al-n}$.
Let $1\le p<n/\al$ and $0\le\la<n$. Then the following linear operators
are bounded:

(1) for $1<p<(n-\la)/\al$,
\[
I_{\al}:M^{p,\la}(\R^{n})\to M^{\frac{(n-\la)p}{n-\la-\al p},\la}(\R^{n}),
\]

(2) for $0<\al<n-\la$,

\[
I_{\al}:M^{1,\la}(\R^{n})\to M_{\ast}^{\frac{n-\la}{n-\la-\al},\la}(\R^{n}).
\]
\end{proposition}

\end{document}